\documentclass[12pt]{amsart}
\usepackage{amsfonts,amssymb,amsmath,amscd,mathtools,amsthm,xcolor,mathrsfs,subfiles,tikz-cd,tikz}
\usepackage[backref,colorlinks, linkcolor=blue, citecolor=magenta]{hyperref}
\usepackage{latexsym}
\usepackage{subfigure}
\theoremstyle{plain}
\newtheorem{theorem}{Theorem}[section]

\newtheorem{lemma}[theorem]{Lemma}

\theoremstyle{remark}

\newtheorem{example}[theorem]{Example}
\theoremstyle{definition}
\newtheorem{definition}{Definition}

\newcommand{\paraB} {\mathcal{P}}
\newcommand{\dasa} {D\'a\v sa  \v Severov\'a}
\usepackage{tikz}

\usepackage{devanagari}
\begin{document}
\title{The Folding Mathematics}
\author{Archana S. Morye}
 \email{asmsm@uohyd.ac.in}
 \address{School of Mathematics and Statistics, University of Hyderabad, Hyderabad-500046, India}
 \thanks{The author is supported by `Institute of Eminence University of Hyderabad' (UoH-IoE-RC5-22-003). }

\maketitle

\begin{abstract}
Origami is the art of paper folding, and it borrows its name from two
Japanese words \emph{ori} and \emph{kami}.  In Japanese, {ori}
means folding, and the paper is called {kami}. While origami is
just a hobby to most, there is a lot more to it. If you fold a square
sheet of paper into any of the traditional origami model (for example
the flapping bird) and unfold it, you can see crease patterns. These
crease patterns tell us that there is a lot of geometry hidden behind
the folds.
		
In this article, we investigate the symbiotic relationship between
mathematics and origami.  The first part of this article explores the
utility of origami in education. We will see how origami could become
an effective way of teaching methods of geometry, mainly because of
its experiential nature.  Complex origami patterns cannot be created
out of thin air.  They usually involve understanding deep mathematical
theories and the ability to apply them to paper folding. In the second
part of the article, we attempt to provide a glimpse of this beautiful
connection between origami and mathematics.
\end{abstract}
\keywords{Origami, geometry, paper folding, fold-able numbers, tree-maker, cubic polynomials.}

\section{Introduction}

Origami is a technique of folding paper into a variety of decorative
or representative forms, such as animals, flowers etc.  The origin of
origami can be traced back to Japan. Originally, the art of paper
folding was called as \emph{orikata}, the craft acquired its current name in
1880 \cite{OrigHist}.

\begin{figure}
\begin{center}
\subfigure[Japanese Origami]{\includegraphics[width=36mm]{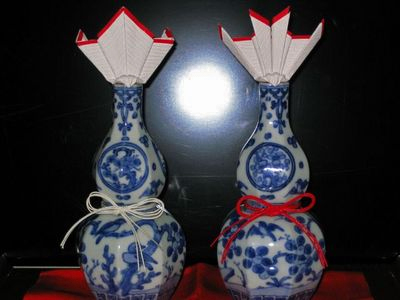}\label{fig:sake}}
\subfigure[European Origami ]{\includegraphics[width=36mm]{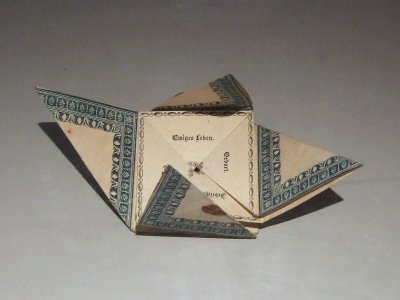}\label{fig:baptist}}
\subfigure[Wet Folding  ]{\includegraphics[width=36mm]{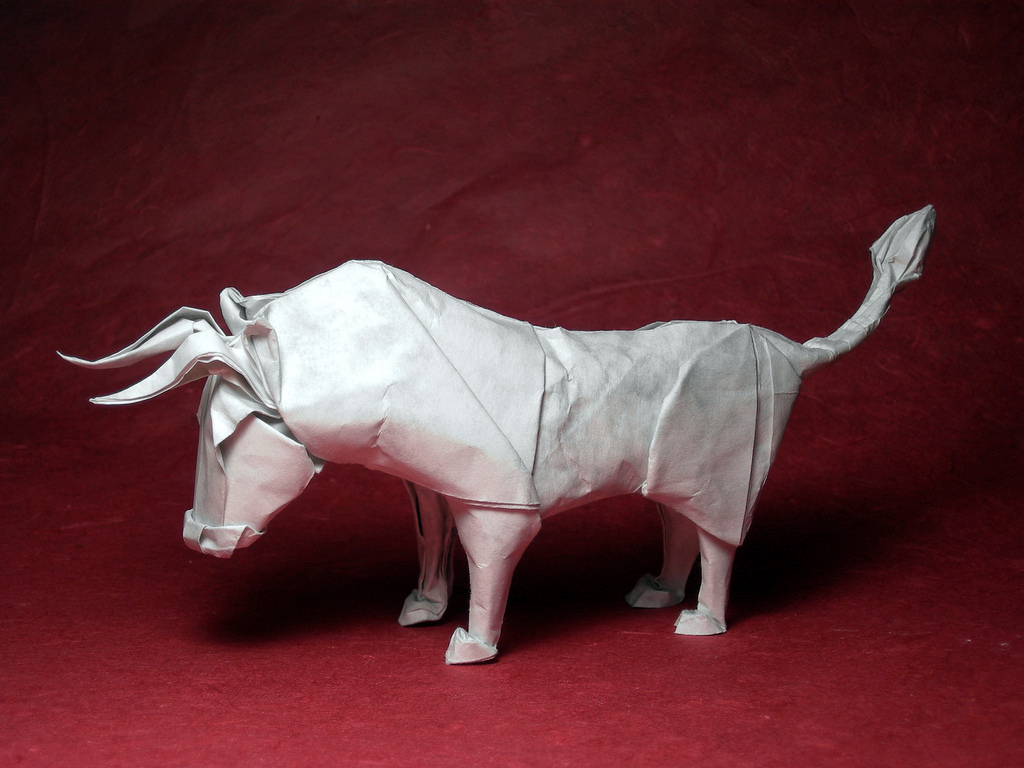}\label{fig:wetfold}}
\caption{}
\end{center}
\end{figure}

The early evidence of origami in Japan suggests that origami was
primarily used as a ceremonial wrapper called the Noshi. \emph{Noshi} is a
wrapper which is attached to a gift, expressing good wishes (similar
to greeting cards of today).  A popular such Noshi is a pair of paper
butterflies known as \emph{Ocho} and \emph{Mecho} that were used to
decorate sake bottles (see Figure \ref{fig:sake}\footnote{Pic source:
\texttt{https://www.origami-resource-center.com/regular-mecho.html}}). 
Origami initially was an art of the elite, mainly because the paper
was a luxury item. As the paper became more accessible, origami also
became a well-practiced art.

The practice of origami can also be traced to Europe, the baptismal
certificates issued during the sixteenth century were folded in a
specific way. (see Figure \ref{fig:baptist}\footnote{Pic source:
 \texttt{https://www.origami-resource-center.com/history-of-origami.html}}).
Here, the four corners of the paper was folded repeatedly to the
center.  Interestingly this techniques is very different from the ones
used in Japan.  It is said that such a crease pattern closely
resembles an old astrological horoscopes.  For this reason, historians
believe that folding in Europe developed more-or-less independently
\cite{OrigHist}.  Some of the popular origami models from Europe are
the Pajarita, the Cocotte and the boat .

 \subsection{The modern origami}
 
The modern era of origami can mainly be attributed to the grandmaster of origami,  Akira Yoshizawa. It was because of his relentless efforts that  origami has transformed into a living art, from  a mere craft. Apart from contributing towards developing more than 50,000 origami models, he also pioneered the wet-folding technique (see Figure \ref{fig:wetfold}\footnote{Pic source: \texttt{https://en.wikipedia.org/wiki/Wet-folding}}). This technique involves slightly dampening the paper before using it for folding. This technique allowed the paper to be manipulated more easily, resulting in models with rounded and  sculpted looks. The famous Yoshizawa-Randlett diagramming system is also his invention. Since the introduction of this system, origami has seen several advancements. Origami is now one of the well-established topics of research in major universities across the world.
 
Origami as it exists today has several variations, we list some of them below.
\begin{enumerate} 
 
 \item \emph{Pure Origami} :
This form of origami is arguably one of the oldest and well studied form,  the models here are made from a single square sheet of paper, without the use of scissors and glue. Coloring the final model is also strictly discouraged.   
 Several models along with instructions can be found in \cite{HappyFold}. My personal favorites and recommendations include the traditional crane, the swallowtail butterfly and the fawn models.
 
 A more stringent variation  is the \emph{Pureland Origami} developed by John Smith. This version disallows certain types of folds allowed by the pure origami.  Interestingly this type of origami is considered to be disable friendly.

 \item \emph{Action Origami}:  This form of origami involves developing models that can be animated. The most famous among the models of this type is the flapping bird. The bird flaps its wings when its tail is moved. Other interesting models of this type include origami airplanes and  the instrumentalist created by Prof. Robert Lang.
\begin{figure}
 	\begin{center}
 	 	 	\subfigure[ Fujimoto Hydrangeas]{\includegraphics[width=31mm]{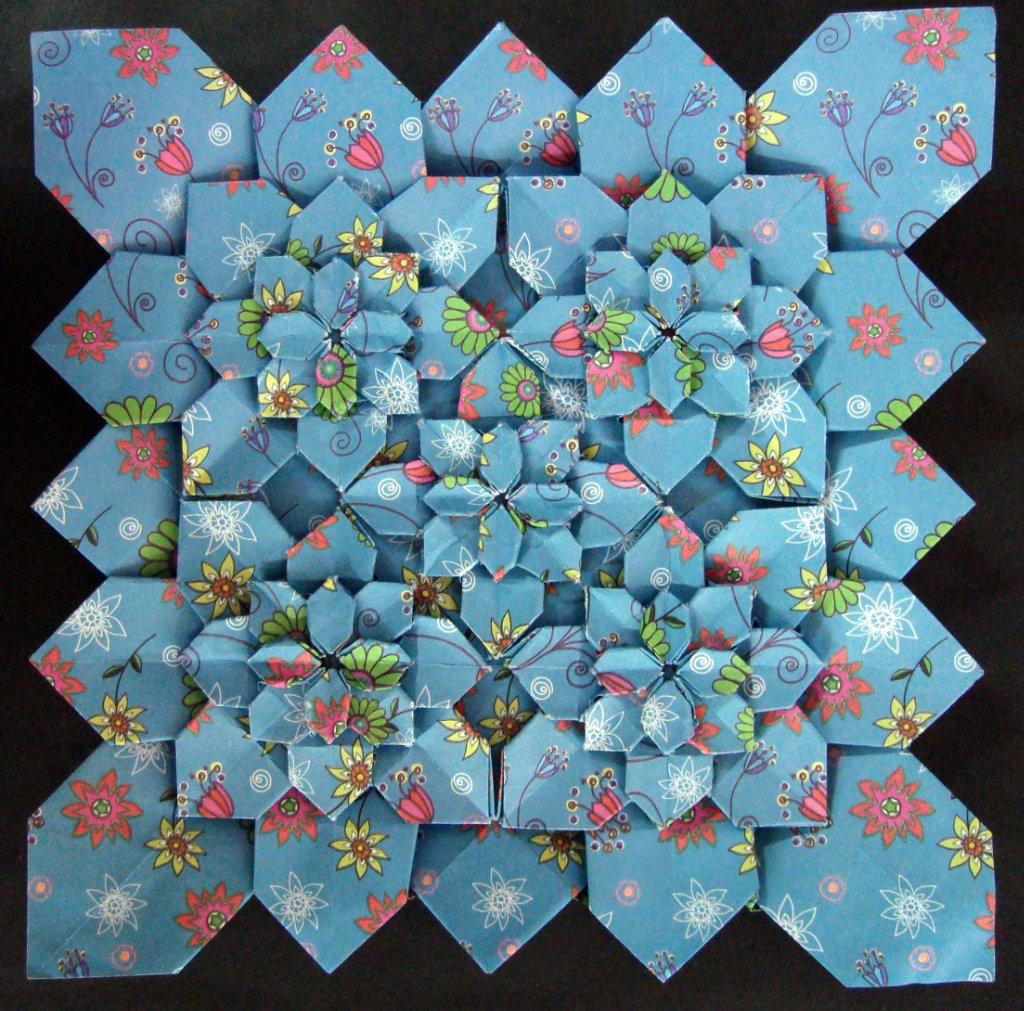}\label{fig:fuji}}
	\subfigure[Kusudama ]{\includegraphics[width=37mm]{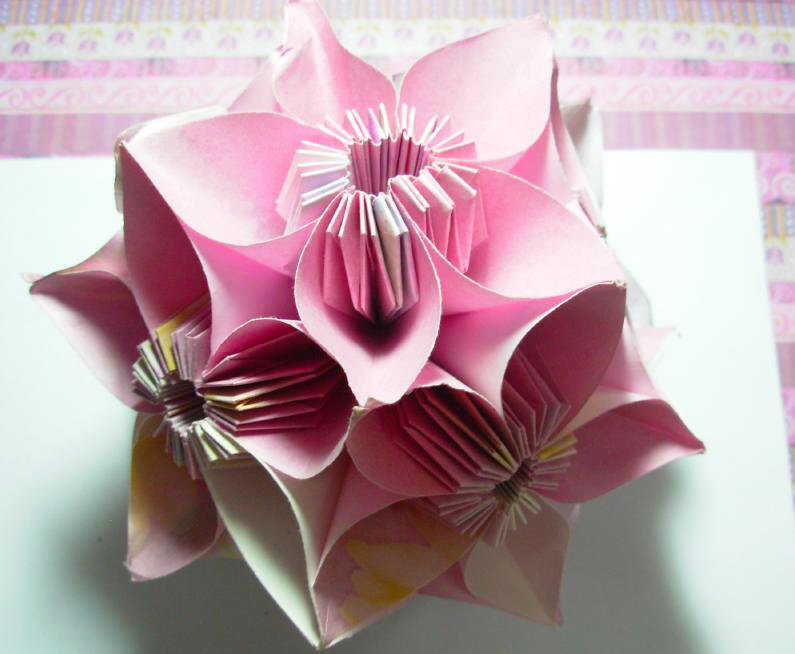}\label{fig:test}}
	\subfigure[knotology  torus ]{\includegraphics[width=41mm]{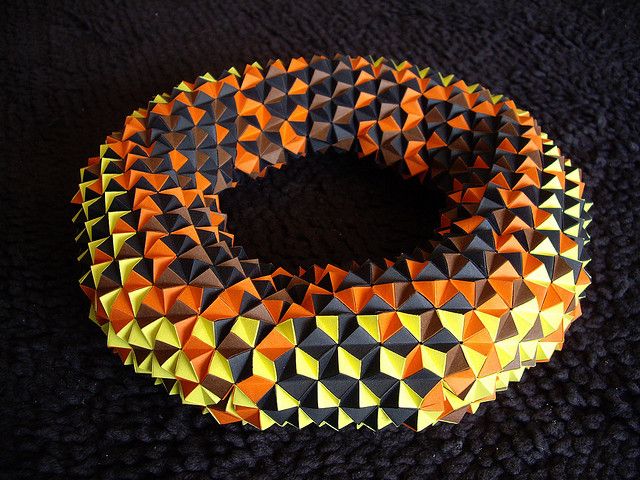}\label{fig:konto}}
\end{center}
 		\caption{}
 \end{figure}

 \item \emph{Modular Origami:}
 In this form of origami, several identical units are folded and then assembled into a more complex origami model. An of this type is the Kusudama flower.  In this model, sixty identical units are folded and arranged into twelve flowers. These flowers are then arranged such that they form a regular dodecahedron,  the pieces are held in place using glue or a thread.  Many models have pocket and flap in each unit, so that units bind together without a glue or a thread.
 
 \item \emph{Origami Tessellation:} 
An origami tessellation is created by repeating a  pattern multiple times in all the directions, and this creates a mosaic. The kind of folds in this process predominantly includes pleats and twists. The invention of this technique can be attributed to Shuzo Fujimoto. This type of origami has an additional feature, and they produce a beautiful effect when they are backlit (when they are held against the light). Fujimoto Hydrangeas (see Figure \ref{fig:fuji}) is an interesting model of this type.
 
\item \emph{Strip Origami}:
 Strip folding is a technique that involves both paper folding and paper weaving.  A fascinating model of this type is that of knotology  torus (see Figure \ref{fig:konto}\footnote{Pic source:\texttt{https://www.flickr.com/photos/dasssa/3426754850}}) by \dasa. 
\end{enumerate}

While origami certainly has evolved as an amazing art, its applicability has also been phenomenal.  Origami-inspired techniques are being sought after by almost every engineering field, ranging from space science to medical equipment and even automobile manufacturers \cite{ScOrigami, IncrOrigami}. In medicine, origami techniques are often applied to stent designs. Stents are collapsible tubes that can be inserted into a patient's veins or arteries. When deployed, the stent expands to open the veins or arteries to improve blood flow. Origami design techniques are instrumental in developing thin and small stents. NASA's James Webb Telescope (JWST), the planned successor of Hubble space telescope, is a rather sizable infrared space telescope with a primary mirror of $6.5$-meters. Origami techniques are being deployed to fold such a large telescope compactly so that it can be airlifted to space, where it can be unfolded again. Automobile manufacturers are pursuing efficient flat-folding techniques for airbags so that it occupies less space and yet unfurls quickly enough when needed.  

The use of origami techniques in the science and engineering field is endless. Hence it is only prudent to study and understand origami as a science. One of the aims of this article is to illustrate the deep connections of origami with mathematics. The article is divided into two parts; the first part deals with using origami as a means to understand mathematics, geometry in particular. For example, a proof for the Pythagoras theorem merely is folding a square sheet of paper. Trisecting a line using just folds is possible. Interestingly, even trisecting an angle can be achieved. The latter is of significant interest due to its impossibility within the realms of Euclidean geometry. The second part of the article briefly dwells upon the need for mathematics to design complex origami models.

Before we delve into the technical details, I would like to share my experience with the folding and origami community.  While this subsection is unconventional,  the reason I write this is to get an opportunity to acknowledge individuals who helped and inspired me during my journey. I also hope that this will nudge others to start their journey into the world of folding.

\subsection{My Origami Journey}

My first exposure to origami was like any other kid when I learned to fold simple models like box, flower, purse, and so on in school. The first complicated model I folded was that of a fish, from a magazine I found in my relatives' place. Learning this model gave me a huge sense of accomplishment and satisfaction. My first formal tryst with origami was during my undergraduate days when I chanced to find a Marathi and English origami series written by Indu Tilak. This series is part of the textbooks prescribed by Maharashtra board for primary education. Even though these are school books, they are quite interesting and extensive. The book includes the necessary foundations to start creating basic and complex structures in origami.

While I was helped and inspired by many in the origami community, one who was extremely helpful and kind to me was \dasa. Papers are the lifeline for origami artists. Some particular models require individual papers. \dasa\ was kind enough to help me obtain certain papers that were difficult to get in India. She has also helped me fold some intricate origami tessellations.

These are some other origami books that one may wish to read.

\begin{itemize}

\item Origami Tessellations: Awe-Inspiring Geometric Designs, by Eric Gjerde.
\item Origami Boxes by Tomoko Fuse.
\item Origami Butterflies by Micheal LaFusse.
\item Origami Journey: Into the Fascinating World of Geometric Origami, by D\'a\v sa  \v Severov\'a.
\item Origami Inspiration by Meenakshi Mukerji.
\end{itemize}
The  website \texttt{www.happyfolding.com}, owned by Sara Adam, is a one-point source of abundant information, beginners or otherwise.

\section{Origami for Mathematics}

Mathematics has unfortunately attained the notoriety of not just dreading the youngsters but also the adults. How many times have we not heard the phrase, ``Thank god, I do not have to study mathematics anymore." This feeling has been captured aptly by the following  Marathi couplet.

\begin{quote}
	{\dn   BolAnAT u@yA aAh\? gEZtAcA pp\?r\\
		poVAt mA\324wyA k\30Fw y\?Un \7{d}K\?l kA r\? Yopr}
\end{quote}

Bholanath is a mystical, mythological, and benevolent Ox. The children pray to him for ill-health, for it is their mathematics paper the next day. It is not very surprising that mathematics is viewed so unfavorably, as it involves many abstract concepts challenging to comprehend. Even a simple definition such as an area or volume is very abstract, for it is a measurement irrespective of the shape. One reason why mathematics lacks the popularity of the other sciences is that it lacks visual appeal. Prof. John J Hopfield (the recipient of the ICTP Dirac medal, 2001)  in one of his articles \cite{Hop2004} wrote that the reason for him being a scientist was because of the encouragement he received to do experiments. Labs and experiments are never associated with mathematics. While one could still argue that mathematics does provide the same experience through puzzles and problems, my personal experiences and memories negate such claims. Solving mathematical problems in the current day scenario has degenerated to learning to apply formulas to score high. Here is where origami can fill in to provide the missing fun.  

\begin{example}
Consider a simple problem of proving that the sum of interior angles of a triangle adds up to $180^\circ$. While there are many ways to prove this, the  folding in Figure \ref{fig:angsum} demonstrates this clearly and crisply.
	\begin{figure}[!htbp]
		\centering
		\begin{tikzpicture}[scale=1.3]
		\thicklines 
		\draw[thick] (0,0)--(3,0)--(1,2)--(0,0);
		\draw [dashed] (.5,0)--(.5,1)--(2,1)--(2,0);
		\node at (-.1,-.1) {\small{$A$}};
		\node at (3.1,-.1){\small{$B$}};
		\node at (1.1,2.1){\small{$C$}};
		\draw [dotted] (4,0)--(7,0)--(5,2)--(4,0);
		\draw [thick] (4.5,0)--(4.5,1)--(6,1)--(6,0);
		\draw [thick] (4.5,1)--(5,0)--(6,1);
		
	\draw[thick](4.5,0)--(6,0);;
		\draw[thick] (0.1,0) arc [radius=.1, start angle=0, end angle=63.45];
		\draw[thick](2.9,0) arc [radius=.1, start angle=180, end angle=135];
		\draw[thick](.97,1.94) arc[radius=.1, start angle=242.45, end angle=300];
		\draw[thick](5.1,0) arc[radius=.1, start angle=0, end angle=180];
		\draw[thick](0.22,0.54)--(0.27,0.46);
		\draw[thick](.72,1.54)--(.78,1.46);
		\draw (2.46,.46)--(2.54,.54);
		\draw (2.41,.51)--(2.49,.59);
		\draw (1.46,1.46)--(1.54,1.54);
		\draw (1.41,1.51)--(1.49,1.59);
		
		\node at (5,-.1) {\tiny{$ A,B,C$}};
		\node at (.2,.1){$\alpha$};
		\node at (1,1.75){$\beta$};
		\node at (2.75,.1){$\gamma$};
		
		\draw [dotted,<->] (1,2)--(1,0);
		\node at (1.1,1.4){$h$};
		\draw [dotted, <->] (0,-0.1)--(3,-0.1);
		\node at (1.2,-0.2){$b$};
		
		\draw [dashed, <-] (4.7,0.2) to[out=120,in=50]  (4.2,0.2) ;
		\draw [dashed, <-] (5.7,0.2) to[out=20,in=110]  (6.5,0.1) ;
		\draw [dashed, <-] (5.1,0.6) to[out=50,in=300]  (5.2,1.6) ;
		
		\node at (4.8,.1){\tiny{$\alpha$}};
		\node at (5,.2){\tiny{$\beta$}};
		\node at (5.2,.1){\tiny{$\gamma$}};
		\node at  (4.3,.5){\tiny{$h/2$}};
		\node at (5,1.1){\tiny{$b/2$}};
		\end{tikzpicture}
		\caption{}\label{fig:angsum}
	\end{figure}
Here the desired rectangle is identified, and the cones of the triangle are folded so that their tips meet. Since the angles $\alpha, \beta, \gamma$ covers a straight line, this proves that their sum is indeed $180^\circ$. Now the same folding also provides us with the argument of why the area of a triangle is $1/2 \times \text{base } \times \text{ height}$. Notice that the fold covers the rectangle with height $h/2$ and length $b/2$ two times. So the area of the triangle is two times the area of this rectangle. This immediately provides us with the relation $$\text{ Area}(\triangle ABC) = 2 \times (b/2) \times (h/2) = 1/2 \times b \times h.$$
\end{example}

This demonstrates that every folding has a deep mathematical connection, possibly many. Discovering them depends on the creativity of the folder. After all, creativity is in the eye of the beholder, beautifully summarized by our beloved Dr.~A.~P.~J. ~Abdul~Kalam as

\begin{quote}

{\bf	``Creativity is seeing the same thing but thinking differently." }
\end{quote}

We will demonstrate the deep connection of origami with mathematics through another example, this time, the famous Pythagorean Theorem. The Pythagorean theorem is one of the oldest known theorems and was studied by Babylonian, Egyptian,  Indian, and Greek mathematicians centuries earlier. It states that the square of the hypotenuse of any right-angled triangle is equal to the sum of the squares of the other two sides. This geometric theorem probably has the most number of proofs. The standard proof which is given in the most high-school books is using similar triangles.  Another exciting proof is as follows. Let $\triangle ABC$ be any right-angle triangle with $c$ as its hypotenuse and its other two sides being $a,b$. Let $C$ be a square-shaped bucket with length and height being $c$-units and its width $1$-unit. Similarly, let $A$ ( respectively $B$) be a square-shaped bucket with length and height $a$ (respectively $b$) and with width exactly $1$-unit. It can be demonstrated that the $C$ bucket can be filled using the water in buckets $A$ and $B$, respectively. While this is undoubtedly a fun proof, the knowledge of volumes is needed to understand the proof. Further conducting such an experiment in a class is cumbersome.
\begin{example}
 Now consider the folding in Figure \ref{fig:pyth}, it clearly demonstrates the proof of Pythagoras Theorem. 
\begin{figure}[!htbp]
	\centering
	\begin{tikzpicture}[scale=.5]
	\thicklines
	\draw[thick] (0,0)--(7,0)--(7,7)--(0,7)--(0,0);
	\draw[dashed] (3,0)--(7,3)--(4,7)--(0,4)--(3,0);
	\node at (-.4,2){$a$};
	\draw [<->](-.2,0)--(-.2,4);
    \node at (-.4,5.5){$b$};
	\draw [<->](-.2,4)--(-.2,7);
\node at (1.9,2.4){$c$};
	\draw[<->](3.2,.15)-- (0.2,4.15); 
	\draw[dotted] (9,0)--(16,0)--(16,7)--(9,7)--(9,0);
	\draw[thick] (12,0)--(16,3)--(13,7)--(9,4)--(12,0);
			\draw [dashed, ->] (1,0.8) to[out=100,in=100]  (3.1,1.2) ;
		\draw [dashed, ->] (5.9,.8) to[out=70,in=90]  (3.8,1.2) ;
		\draw [dashed, ->] (5.6,5.8) to[out=-40,in=-20]  (4.5,4.6) ;
		\draw [dashed, ->] (1.2,5.6) to[out=-100,in=-150]  (2.6,4.5) ;
			\end{tikzpicture}
	\caption{} \label{fig:pyth}
\end{figure}
In the figure, we take a square sheet of paper and mark out four right-angled triangles of equal sizes along the four corners. Each of these triangles has $c$ as its hypotenuse and $a,b$ as the size of its other sides. Each of these triangles has its hypotenuse in the inner part of the square, touching each other. Dotted lines in the figure denote these. Folding along the hypotenuse gives us a square of length $c$. Notice that the area of such a square is $c^2$, the area of the original square we started with was $(a+b)^2$. What was folded in were $4$ triangles each of area $1/2 \times a \times b$. With this, we get the following equation, which also proves the Pythagoras Theorem.

$$ c^2 = (a+b)^2 - (4 \times ( 1/2 \times a \times b)) = a^2 + b^2. $$

\end{example}

While we could go on demonstrating the utility of origami in proving theorems involving simple properties, one may ask whether origami can also be used to solve more involved problems and theorems. For this, we need to formalize folding, i.e., make precise what kind of folds are allowed and what are not. This would allow us to investigate the constructible geometric objects through origami. This is in the same lines as the classical \emph{ruler-and-compass} construction.

\subsection{Huzita-Hatori Axioms}
 
The formal axioms for origami is given by the \emph{Huzita–Hatori axioms}. We briefly recall them here and direct the interested readers to \cite{L} for a comprehensive coverage on the subject.
 
 \begin{enumerate}
 	
 	\item Given two distinct points $p_1$ and $p_2$, there is a unique fold that passes through both of them.
 	\item Given two distinct points $p_1$ and $p_2$, there is a unique fold that places $p_1$ onto $p_2$.
 	\item Given two lines $l_1$ and $l_2$, there is a fold that places $l_1$ onto $l_2$.
 	\item Given a point $p_1$ and a line $l_1$, there is a unique fold perpendicular to $l_1$ that passes through point $p_1$.
 	\item Given two points $p_1$ and $p_2$ and a line $l_1$, there is a fold that places $p_1$ onto $l_1$ and passes through $p_2$.
 	\item Given two points $p_1$ and $p_2$ and two lines $l_1$ and $l_2$, there is a fold that places $p_1$ onto $l_1$ and $p_2$ onto $l_2$. 
\item Given a point $p$, and two lines $l_1$ and $l_2$, we can make a fold perpendicular to $l_2$ that places $p$ onto line $l_1$.
 	 \end{enumerate}

\begin{figure}[!htbp]
\begin{center}
\begin{tikzpicture}[scale=.7]
\thicklines
\draw[thick] (0,0)--(2,0)--(2,2)--(0,2)--(0,0);
\draw[fill=black](.5,1) circle(.05);
\node at (.5,.8) {\tiny{$p_1$}};
\draw[fill=black](1.5,1) circle(.05);
\node at(1.5,.8) {\tiny{$p_2$}};
\draw[dashed, thick] (0,1)--(2,1);
\node at (.6,2.2) {\tiny{Axiom 1}};

\draw[thick] (2.5,0)--(4.5,0)--(4.5,2)--(2.5,2)--(2.5,0);
\draw[fill=black](3,1.25) circle(.05);
\node at (3,1) {\tiny{$p_1$}};
\draw[fill=black](4,.75) circle(.05);
\node at(4,.6) {\tiny{$p_2$}};
\draw[dashed, thick] (3,0)--(4,2);
\node at (3.1,2.2){\tiny{Axiom 2}};

\draw[thick] (5,0)--(7,0)--(7,2)--(5,2)--(5,0);
\node at (5.5,1) {\tiny{$l_1$}};
\draw[thick](7,.5)--(5,1.5);
\draw[dashed,thick] (5,0.666)--(7,1.334);
\node at(6,.4) {\tiny{$l_2$}};
\draw[thick] (5.5,0)--(6.5,2);
\node at (5.6,2.2){\tiny{Axiom 3}};

\draw[thick] (7.5,0)--(9.5,0)--(9.5,2)--(7.5,2)--(7.5,0);
\node at (9,.8) {\tiny{$l$}};
\draw[thick](7.5,1)--(9.5,1);
\node at (8.3,1.5) {\tiny{$p$}};
\draw[fill=black](8.5,1.5) circle(.05);
\draw[dashed,thick] (8.5,2)--(8.5,0);
\draw[thick] (8.5,1.2)--(8.7,1.2)--(8.7,1);
\node at (8.1,2.2){\tiny{Axiom 4}};

\draw[thick] (10,0)--(12,0)--(12,2)--(10,2)--(10,0);
\draw[thick](10,.5)--(12,.5);
\node at (11.5,.7) {\tiny{$p_1$}};
\draw[fill=black](11.5,1) circle(.05);
\node at (10.75,1.5) {\tiny{$p_2$}};
\draw[fill=black](10.5,1.5 ) circle(.05);
\draw[->] (11.5,1) arc(90:180:.5);
\draw[dashed,thick] (10,2)--(12,0);
\node at (10.5,.3) {\tiny{$l_1$}};
\node at (10.6,2.2){\tiny{Axiom 5}};

\draw[thick] (12.5,0)--(14.5,0)--(14.5,2)--(12.5,2)--(12.5,0);
\draw[thick](12.5,1)--(13.5,.5)--(14.5,1);
\node at (14,1.5) {\tiny{$p_1$}};
\draw[fill=black](14,1.2) circle(.05);
\node at (13.25,1.5) {\tiny{$p_2$}};
\draw[fill=black](13,1.5 ) circle(.05);
\draw[->] (13,1.5) to [out=0, in=90] (13,.75);
\draw[->] (14,1.2) to [out=0, in=90] (14,.75);
\draw[dashed,thick] (12.5,.75)--(14.5,.75);
\node at (13.1,2.2){\tiny{Axiom 6}};

\draw[thick] (15,0)--(17,0)--(17,2)--(15,2)--(15,0);
\node at (16.5,.4) {\tiny{$l_2$}};
\draw[thick](15,.25)--(17,.25);
\node at (16.5,1.6) {\tiny{$l_1$}};
\draw[thick](15.75,2)--(17,.75);
\node at (15.7,1.4) {\tiny{$p$}};
\draw[fill=black](15.5,1.25) circle(.05);
\draw[fill=black](16.5,1.25) circle(.05);
\draw[dashed,thick] (16,2)--(16,0);
\draw[thick] (16,.45)--(16.2,.45)--(16.2,.25);
\node at (15.6,2.2){\tiny{Axiom 7}};
\end{tikzpicture}
\caption{}\label{fig:axioms}
\end{center}
\end{figure}

The relevant question here is, what can these postulates achieve can they achieve anything more in comparison to the classical ruler-and-compass.
While the first four postulates are self-explanatory, we will examine the fifth and sixth postulates closely. Interestingly the fifth postulate can be used to solve a quadratic equation
and the sixth a cubic equation. We will demonstrate the latter in the sequel. We first prove the following lemma which states that the dotted line in the Axiom $5$ of Figure \ref{fig:axioms}
is actually a tangent to a parabola with its focus on $p_1$ and the directrix on $l_1$. Recall that a parabola is those set of points (called the \emph{locus}) that are equidistant from a fixed point (called the \emph{focus}) and a fixed-line (called the \emph{directrix}).

\begin{lemma}\label{lem:parabola}
	Given two points $p_1$ and $p_2$ and a line $l_1$, the fold  $\ell$ that places $p_1$ onto $l_1$ and passes through $p_2$ is tangent to the parabola $\mathcal{P}$ defined by the focus $p_1$ and the directrix $l_1$.
\end{lemma}

\begin{proof}
	To prove the lemma, we prove that there is a unique point $x$ on the line $\ell$ which is equidistant from $p_1$ and $l_1$. 	By definition, this point will lie on the parabola. Since this point is unique,  no other points of the line will lie on the parabola. Hence, the line will be tangential to $\mathcal{P}$. 
	
	For this, we prove two things. Firstly we prove that there is a point $x$ in $\ell$, which is equidistant from $p_1$ and $l_1$. Here the distance of $x$ from $l_1$ is given by the shortest one. Then we will prove that for any point $y$ in $\ell$ such that $y \neq x$, its shortest distance from $l_1$ is not equal to its distance from $p_1$.
	\begin{figure}[!htbp]

\begin{minipage}{.35\textwidth}
		
\begin{tikzpicture}[scale=.95]
\thicklines
	\draw[thick] (0,0)--(4,0)--(4,4)--(0,4)--(0,0);
\draw(1,0)--(1,2)--(3,2);
\draw[dashed,thick](3,0)--(0,3);
\node at (.1,3.2){\tiny{$\ell$}};
\node at (-.2,3) {\tiny{$p_2$}};
\node at (3.5,-.3){\tiny{$l_1$}};
\draw[fill=black](3,2) circle(.05);
\node at (3.3,2) {\tiny{$p_1$}};
\draw[fill=black](1,0) circle(.05);
\node at (1,-.25){\tiny{$\bar{p_1}$}};
\draw(1,0)--(3,2);
\draw[fill=black](2,1) circle(.05);
\node at (2,1.3){\tiny{$o$}};
\draw[fill=black](1,2) circle(.05);
\node at (1,2.2) {\tiny{$x$}};
\draw (.5,0)--(.5,2.5)--(1,0);
\draw (1.9,1.1)--(2,1.2)--(2.1,1.1);
\draw[fill=black](.5,2.5) circle(.05);
\node at(.6,2.6){\tiny{$y$}};
\draw[fill=black](.5,0) circle(.05);
\node at(.5,-.25){\tiny{$p_y$}};
\draw(.5,.1)--(.6,.1)--(.6,0);
\end{tikzpicture}
\caption{}\label{fig:lemma1}
\end{minipage}
\qquad
\qquad
\begin{minipage}{.45\textwidth}
	\begin{tikzpicture}[scale=.45]
	\thicklines
	\draw[dotted] (-4,0)--(6,0);
	\draw[thick](-4,-1)--(6,-1);
	\draw[dotted](0,4)--(0,-6);
	\draw[thick](2,4)--(2,-6);
	\node at (2,4.3){\tiny{$l_2:x=2$}};
	\node at (7.5,-1.5) {\tiny{$l_1:y=-1$}};
	\node at (0.5,1.3){\tiny{$p_1=(0,1)$}};
	\draw[fill=black](0,1) circle(.05);
	\node at (5.4,-.6) {\tiny{$\bar{p_1}=(4,-1)$}};
	\draw[fill=black](4,-1) circle(.05);
	\draw[thick](0,1)--(4,-1);
	\draw[thick,dashed](4,4)--(-1,-6);
	\node at (5.7,4.2){\tiny{$\ell:y=2x-4$}};
	\node at (-2.4,-2.7){ \tiny{$p_2=(-2,-3)$}};
	\draw[fill=black](-2,-3) circle(.05);
	\draw[thick] (2,-5)--(-2,-3);
	\draw[fill=black](2,-5) circle(.05);
	\node at(3.8,-5)  {\tiny{$\bar{p_2=(2,-5)}$}};
	
	\end{tikzpicture}
	\caption{$x^3+0 x^2-3x-2$}\label{fig:lemma2}
\end{minipage}

\end{figure}
	Without loss of generality, we will assume that the line $l_1$ is the bottom edge of a square and that the point $p_2$ lies on the left edge of the square. Observe that given any $l_1,p_1, p_2$, one could arrange them in this manner in an appropriate large enough square. Recall that the axiom guarantees that the line $\ell$ folds $l_1$ in such a way that a point in it coincides with $p_1$. Let this point on $l_1$ be $\bar{p_1}$. Now draw a line perpendicular to the line $l_1$, starting at $\bar{p_1}$, let this line intersect $\ell$ at $x$. We claim that this intersection point $x$ is the required point. That is, it is equidistant from both $l_1$ and $p_1$. This is easy to observe since $\bar{p_1}$ coincides with $p_1$ when the paper is folded along $\ell$. This immediately implies that $\overline{xp_1} = \overline{x \bar{p_1}}$.
	
	For proving the second part, we again use the fact that $\bar{p_1}$ coincides with $p_1$ when the paper is folded along $\ell$. This immediately also tells us that for any point $y$ on the line $\ell$, $\overline{y p_1} = \overline{y \bar{p_1}}$. Notice that the shortest distance from any point $y$ in $\ell$ to the line $l_1$ is obtained by dropping a line perpendicular to $l_1$ from $y$. Let this line intersect $l_1$ at $p_y$. Now for any point $y \neq x$, consider the right angle triangle $\triangle y p_y \bar{p_1} $. In this triangle since $\overline{y \bar{p_1} }$ is the hypotenuse, it follows that  $\overline{y \bar{p_1} } > \overline{y {p_y} }$ 	
\end{proof}

Notice that in the proof above, the role of $p_2$ is non-existent. The fact that $\ell$ is a tangent to the parabola $\mathcal{P}$ is invariant to $p_2$. Also, notice that there are infinite tangents to the parabola, one for every point on it. This collection of tangents is called the \emph{tangent bundle}. The point $p_2$ determines a unique tangent from such a tangent bundle.

Using the same technique of Lemma \ref{lem:parabola}, one could also prove that the fold $\ell$ obtained in Axiom $6$ is a line tangent to two parabolas $\mathcal{P}_1, \mathcal{P}_2$, determined by the focus and directrix $p_1$, $l_1$ and $p_2$, $l_2$ respectively.
We will instead illustrate how to use Axiom $6$ to solve a cubic equation. More specifically, we will show how to obtain a solution for equations of the form $x^3 + a x^2 + b x +c = 0 $ where $a$, $b$, $c \in \mathbb{R}$.
The idea here is to use the points $p_1 =  (a,1)$, $p_2 = (c,b) $, and the lines $l_1$ given by $y = -1 $ and $l_2$ given by $ x = -c$ in Axiom $6$ and show that the slope of the fold obtained is a solution to the equation.

\begin{theorem}
Let $x^3 + a x^2 + b x +c = 0 $  be any cubic equation with $a,b,c \in \mathbb{R}$. Consider a large enough square paper with origin at its center, let $p_1$ and $p_2$ be two point in it with the coordinates given by $(a,1)$ and $(c,b)$ respectively. Let $l_1$ and $l_2$ be lines defined by the equation  $y = -1 $ and $ x = -c$ respectively. Then the slope $t$ of the line $\ell$ that folds the point $p_1$ onto $l_1$ and $p_2$ onto $l_2$ is a solution to the given cubic equation.
\end{theorem}

\begin{proof}
We wish to prove that any solution to the slope of the line $\ell$ obtained as a result of applying Axiom $6$ to the points and lines specified below, is a solution to the given cubic equation $x^3 + a x^2 + b x +c = 0 $. The points are specified as  $p_1 =  (a,1)$ and $p_2 = (c,b)$ and the lines given by the equation $l_1 : y = -1 $ and $l_2: x = -c$. Let the equation defining the line $\ell$ be $\ell : y = tx + u$. Recall that, $t$ is the slope and $u$ is the $y$ intercept here.
	 
Firstly notice that in both Axiom $5$ and $6$, the point $p_1$ is placed onto the line $l_1$. By Lemma \ref{lem:parabola},  and the fact that the lemma is invariant the point $p_2$ in Axiom $5$, we obtain that the line $\ell$ is tangent to the parabola defined by its focus on $p_1$ and the directrix on the line $l_1$. Let this parabola be $\paraB$. Using the focus given by  $p_1$ and the directix given by line $l_1$, we obtain the equation of $\paraB$ as 
		\begin{equation} \label{eq1}
 		y = \frac{1}{4} (x-a)^2.
 		\end{equation}
		Since the line $\ell$ is tangent to the parabola, we obtain the equation of its slope as $t = \frac{\partial y}{\partial x} = \frac{1}{2} (x-a)$. Let $(x_1,y_1)$ be the point on $\ell$ where it makes contact with the parabola. Evaluating the equation of slope at this point provides us with:
	\begin{equation}\label{eq11}
	  t =  \frac{1}{2}(x_1 - a).
	  \end{equation}
	 	 Since we know the slope of $\ell$ and also that $(x_1,y_1)$ lies on $\ell$, we obtain the equation of $\ell$ as follows.
	  	\begin{equation} \label{eq2}
	 y = tx -tx_1 + y_1.
	 \end{equation}
	  From this we obtain that $u =  -tx_1 + y_1$. Further from Equation \ref{eq1}, we know that $y_1 = t^2$. We also know the equation of $x_1$ in terms of $t$ ( $x_1 = 2t +a$ from the Equation \ref{eq11}),  as a result we get the following equation.
	  	\begin{equation} \label{eq2}
			 u = - t^2 -ta.
	 \end{equation}
	  Notice that Axiom $6$ mandates that the line $\ell$ is tangential to another parabola $\paraB'$, defined by the focus $p_2 = (c,b)$ and the line and $l_2: x = -c$. The equation of this parabola is given by   
	 $$ x = \frac{1} {4c} (y-b)^2.$$ 
	 Applying the technique seen earlier, we obtain that the equation of the slope of $\ell$ in this case to be $ t = \frac{2c} { (y-b)}$. Using this, we get $ u =  b + \frac{c}{t}$. Equating the value of $u$ obtained in Equation \ref{eq2} with this, we obtain the following equation. 
	 \begin{equation}
	   b + \frac{c}{t} = - t^2 -ta.
	 \end{equation}
	 Notice that this translates to the equation $t^3 + a t^2	 +bt + c =0$, matching with the original equation that we started with, this completes the  proof.	 
	\end{proof}

\section{ Mathematics for Origami}
While in the previous section, we saw that origami could be very handy to learn mathematical concepts, in this section effectively, we will see that the relationship is symbiotic. Many of the origami models require deep mathematical insights and techniques, and we will survey some of them.

\begin{figure}[h!]
	\centering
	\subfigure{\includegraphics[width=40mm]{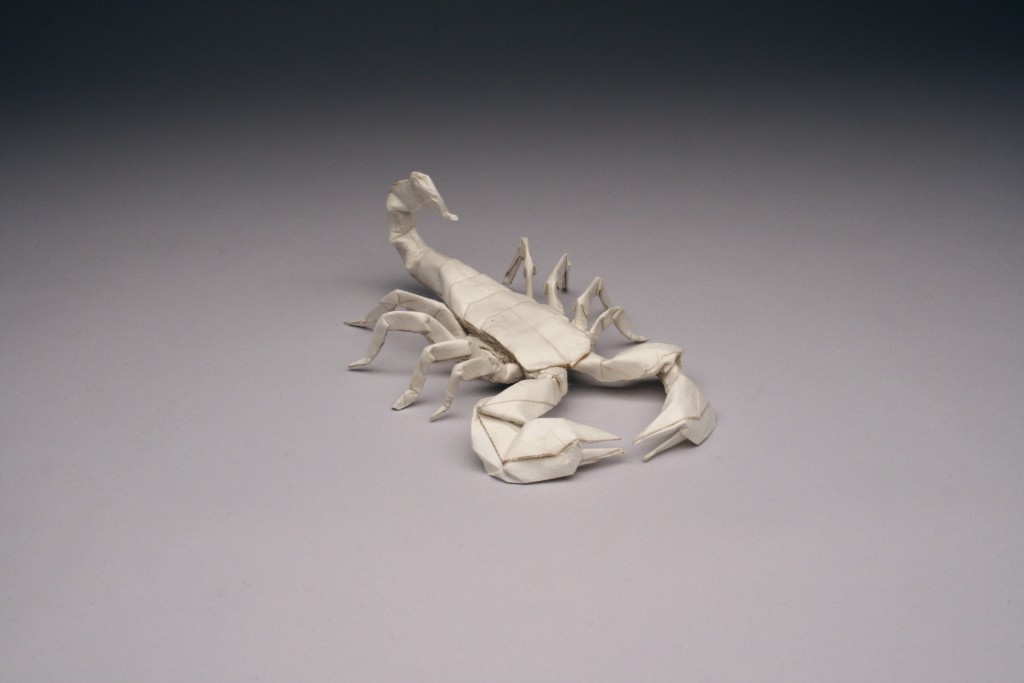}}
	\subfigure{\includegraphics[width=40mm]{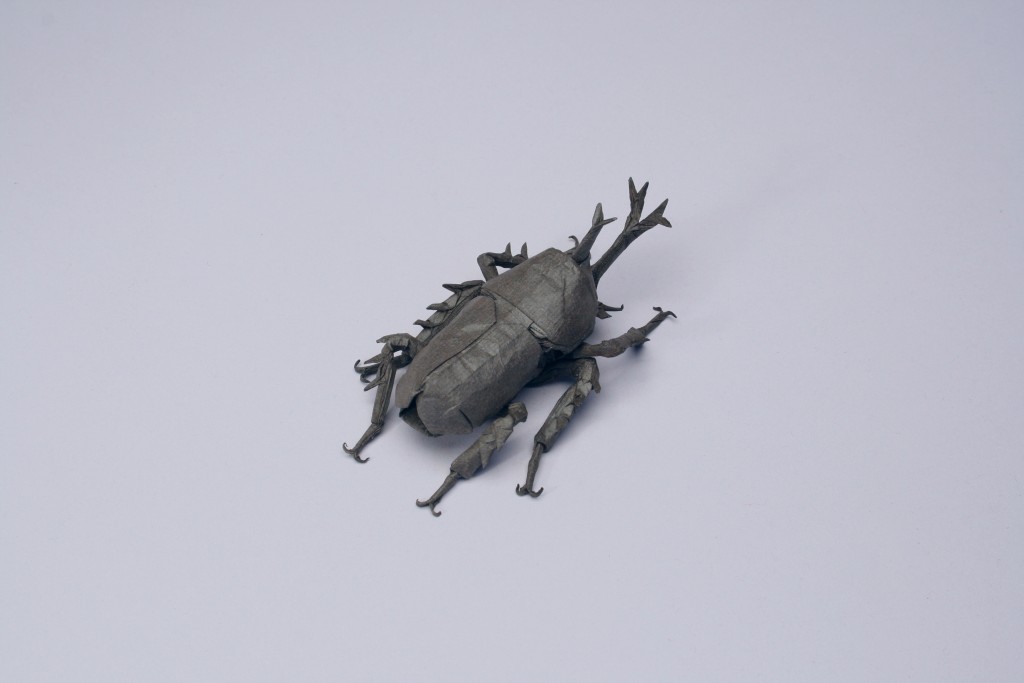}}
	
	\caption{}\label{fig:TreeM}
\end{figure}

The complexity of an origami model depends on the number of attachments it has. For example, an origami model resembling a bird with a head, a tail, and two wings is less complex than a model resembling a beetle with a head, two horns, and six legs. One approach to building an origami model is to come up with what is called a base. A \emph{base} is a geometric object which roughly resembles the end product. Historically most of the origami models were constructed by trial and error basis.  That is, the paper is folded roughly in the direction of the base until the desired objective is achieved. One of the most difficult and important tasks in coming up with an origami construction is to identify how and when to create a fold. This translates to identifying the crease pattern required for the construction. While for simpler models, the trial and error approach is possible, it is highly inefficient for complex models. One important and pertinent question in this regard is whether one can come up with the crease patterns that can be folded into the desired base. The computer programs TreeMaker and Origamizer achieve this objective. Here, we survey the techniques used in the Tree Maker program. We closely adopt the definitions and language used in the TreeMaker manual \cite{TreeMan}. To keep the article simple, we are making our exposition is very brief and informal.  We direct the interested readers to \cite{} for an extensive exposition on the subject. TreeMaker is a computer program developed by a famous origami artist and scientist Robert Lang. The software has been used to develop several complex origami models, see Figure \ref{fig:TreeM}\footnote{Pic source: {\texttt{https://langorigami.com/artworks/}}}.

\begin{figure}
	\centering
		\subfigure{\includegraphics[width=35mm]{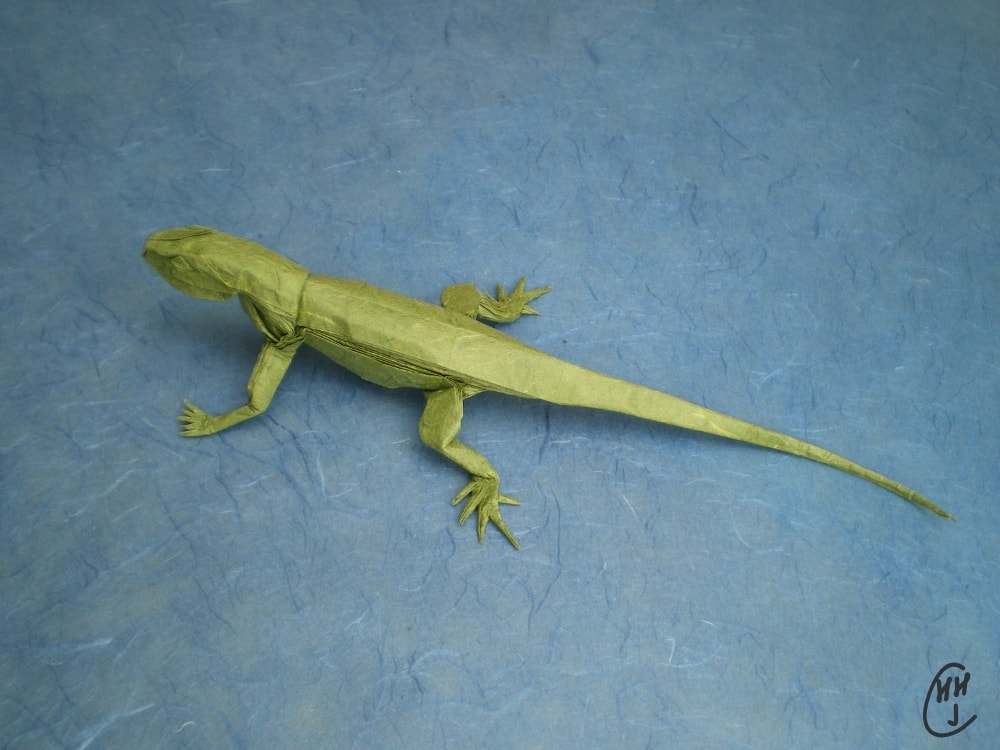}}
		\hspace{75pt}
\begin{tikzpicture}[scale=0.5]
\thicklines 

\draw[thick] (0,0)--(7,0)--(7,5)--(0,5)--(0,0);

\draw[thick] (6,0.5)--(0.75,4.5);

{\fill[black] (6,0.5) circle (4pt);}
\draw (5.3,0.5) node {{\tiny tail}};

{\fill[black] (0.75,4.5) circle (4pt);}
\draw (1.5,4.7) node {{\tiny head}};
\draw (1.6,4.2) node {{\tiny $1$}};

{\fill[black] (1.7,3.55) rectangle (2.0,3.8);}

{\fill[black] (3.5,2.20) rectangle (3.8,2.45);}

\draw[thick] (1.7,2.20)--(1.75,3.7) -- (3.5,3.7);

\draw (2.7,3.9) node {{\tiny $1$}};
\draw (1.5,2.8) node {{\tiny $1$}};

\draw (2.5,2.8) node {{\tiny $2$}};

{\fill[black] (1.7,2.20) circle (4pt);}
{\fill[black] (3.5,3.7) circle (4pt);}

\draw[thick] (3.55,0.90)--(3.55,2.32) -- (5.25,2.32);

\draw (3.4,1.6) node {{\tiny $1$}};
\draw (4.6,2.55) node {{\tiny $1$}};

\draw (4.65,1.2) node {{\tiny $3$}};

{\fill[black] (3.55,0.9) circle (4pt);}
{\fill[black] (5.25,2.32) circle (4pt);}



\end{tikzpicture}	
\caption{} \label{fig:treefig}
\end{figure}

Given any base, let its tree diagram be the graph obtained by shrinking its skeleton to straight lines. In reality, the edges of such a tree diagram can have weights; these weights correspond to the size of foldings required in the base. For example, Figure \ref{fig:treefig}\footnote{Pic source: \texttt{https://origami.me/lizards/}} is the tree diagram of the lizard model next to it. Notice that it has a head, a tail, two forelegs, and two hind legs, not all have the same size. A tree diagram is said to be uni-axial if firstly it is tree-like (i.e., it is connected and has no cycles) and further, it has a main stem, and every branch only originates from this stem. The TreeMaker algorithm works when the required base is of type uni-axial.

The tree diagram is crucial to obtaining the desired algorithm. Every vertex (in a graph sense) of such a tree diagram, which has no outgoing edges, is called the terminal node. In the figure, these are represented by circular nodes. Every other vertex is called the \emph{internal node}; in the figure, these are the rectangular nodes. An important mathematical property that allows a given uni-axial tree diagram to be folded to its base is as follows:\begin{quote}\emph{ If one can find in a square, a set of points, each corresponding to a vertex in the tree diagram such that the distance between any two points is greater or equal to the distance between the corresponding vertices in the tree diagram, then it is possible to fold such a set of points into the required base.}\end{quote}

Notice that the property is existential and does not immediately provide a way to recover the folding pattern. To obtain the crease pattern, the first observation is as follows: \begin{quote}\emph{If the distance between any two terminal points is equal to the distance between their corresponding vertices in the tree diagram, then the line between them is definitely a crease in the final base.}\end{quote} The algorithm crucially identifies points in the square so that maximal such creases occur, partitioning the square into polygons. This already is a computationally hard problem and has deep connections to a famous graph theoretical problem called the \emph{cycle packing problem}.

The algorithm then depends on yet another mathematical property that: creases for each polygon partition of the square can be identified separately. The handling of a triangle and quadrilateral is well known, i.e., there are well-known methods to find creases for such simple polygons. However, other complex polygons can create complexities. One way to get around this is to simplify these polygons by adding extra terminal vertices to the tree diagram. This roughly translates to refining the creases and hence, partitioning of the complex polygons into simpler polygons. The algorithm can now be summarized as follows:

\begin{itemize}
	\item Obtain the tree diagram of the desired base.
	\item Obtain points on the square such that the points have appropriate distances.
	\item Mark all the  creases with minimum lengths.
	\item Identify the polygons partitioning the square, based on the crease marking obtained earlier.
	\item Process the crease pattern for the polygons individually. Simplify the polygons if needed.
\end{itemize}

While this algorithm provides the necessary crease pattern to obtain the base, identifying whether a crease is an inside fold (valley fold) or an outside fold (mountain fold) is another problem entirely. While there are effective heuristics for this, the problem, in general, is still  open.

Another interesting problem that one encounters in origami is that of identifying whether a given crease pattern can be \emph{flat-folded}. An origami model is said to be \emph{flat-folded} if it can be compressed without making any additional creases. The question is whether it can be determined automatically, just by looking at the crease pattern, if it can be flat-folded. While there are results for sub-classes of the crease pattern due to Maekawa Jun and Kawasaki Toshikazu, the general problem still remains open. Thus, there are very many mathematical properties of origami that are being actively studied, some of which are still waiting to be solved. This demonstrates the profound impact that mathematics has on the folding art.

\end{document}